\documentclass{amsart}

\usepackage{amssymb}

\theoremstyle{plain}

\numberwithin{equation}{section}

\newtheorem{theo}[equation]{Theorem}

\newtheorem*{teo}{Theorem A}

\newtheorem{induction}[equation]{}

\newtheorem{hyp}[equation]{Hypotheses}

\newtheorem{lem}[equation]{Lemma}

\newtheorem{prop}[equation]{Proposition}

\newtheorem{claim}[equation]{Claim}
\newtheorem{subclaim}[equation]{Step}

\newcommand{\Irr}{\operatorname{Irr}}
\newcommand{\Lin}{\operatorname{Lin}}

\newcommand{\modu}{\operatorname{mod}}

\newcommand{\Ker}{\operatorname{Ker}}

\newcommand{\core}{\operatorname{core}}

\theoremstyle{definition}

\newtheorem{com}[equation]{}
\begin{document}	
\title{Induction of characters and finite $p$-groups}

\author{Edith Adan-Bante}

\address{University of Southern Mississippi Gulf Coast, 730 East Beach Boulevard,
 Long Beach MS 39560}

\email{Edith.Bante@usm.edu}


\subjclass{20c15}

\date{2005}
\begin{abstract} Let $G$ be a finite $p$-group, where $p$ is an odd prime number,
 $H$ be a subgroup of $G$ and 
$\theta\in \Irr(H)$ be an irreducible character of $H$. Assume also that $|G:H|=p^2$.
Then the character $\theta^G$ of $\,G$ induced by $\theta$ is either a multiple of 
an irreducible character of $G$, 
or has at least $\frac{p+1}{2}$ distinct irreducible
constituents. 
\end{abstract}
\maketitle

\begin{section}{Introduction}

Let $G$ be a finite group. Denote by $\Irr(G)$ the set of 
irreducible complex characters of $G$. 
Through this work,
we use the notation of \cite{isaacs}. In addition, we are going to 
denote by $\Lin(G)=\{\lambda\in \Irr(G) \mid \lambda(1)=1\}$ the set of
linear characters. 

Let $\Gamma$ be a character of $G$. Then $\Gamma$ can be express as a nontrivial
 integral linear combination of distinct irreducible characters of $G$. Denote by
 $\eta(\Gamma)$ the number of distinct irreducible constituents of $\Gamma$.

Let $G$ a finite $p$-group, where $p$ is a prime number,
$H$ be a subgroup of $G$ and $\theta \in \Irr(H)$.
Denote by  $\theta^G$ the character of $G$ induce by $\theta$. If 
$H$ is a normal subgroup, then either
 $\eta(\theta^G)=1$, i.e. $\theta^G$ is a multiple
of an irreducible, or $\eta(\theta^G)\geq p$, i.e. $\theta^G$ is an integral
linear combination of at least $p\,$ distinct irreducible characters of $G$
 (see Lemma \ref{normal}). In
Theorem \ref{extensiondade}, it is shown that given any prime $p>2$ and any integer $l\geq 2$,
 there exist a $p$-group $G$, a subgroup
 $H$ of $G$ with
$|G:H|=p^l$
and $\theta\in \Irr(H)$ such that $\eta(\theta^G)=\frac{p+1}{2}$.
 Therefore Lemma \ref{normal}
does not remain true without the  hypothesis that $H$ is normal in $G$. But given any prime $p>2$ and any integer $n>0$, do there exist a $p$-group $G$, a subgroup $H$ of $G$ and 
$\theta\in \Irr(H)$ with $\eta(\theta^G)=n$? If we also required, in addition, that
 $|G:H|=p^2$ and $1<n<\frac{p+1}{2}$, then the answer is no. More specifically

\begin{teo}
Let $G$ be a finite $p$-group, where $p$ is an odd prime number,
 $H$ be a subgroup of $G$ and $\theta\in\Irr(H)$.
Assume also that $|G:H|=p^2$.
Then either $\eta(\theta^G)=1$ or $\eta(\theta^G)\geq \frac{p+1}{2}$.
\end{teo}  
 
 For a fix prime $p>3$, Theorem A implies that there exists a ``gap" among the possible
 values that $\eta(\theta^G)$ can take for any finite $p$-group $G$, 
 any subgroup $H$ of $G$ with $|G:H|=p^2$, and  any character
 $\theta\in \Irr(H)$.
 But, do there exist a $p$-group $G$, a subgroup $H$ of $G$ and 
$\theta\in \Irr(H)$ with $1<\eta(\theta^G)<\frac{p+1}{2}$ and $|G:H|>p^2$? The
answer is yes. In Theorem \ref{examples2}, given any prime $p$ such that
3 divides $p-1$, we provide a  $p$-group $G$, a subgroup $H$ of $G$ with 
$|G:H|=p^3$ and a character $\lambda\in \Lin(H)$ such that
$\eta(\lambda^G)=\frac{p+2}{3}$.  Does it mean then that, for a fixed prime $p>5$,
 there are no ``gaps"
among the possible values that $\eta(\theta^G)$ can take for any finite $p$-group $G$, 
 any subgroup $H$ of $G$ with $|G:H|=p^3$, and  any character
 $\theta\in \Irr(H)$?
 We do not know
 the answer of that question.
\end{section}
\begin{section}{Preliminaries}

\begin{lem}\label{inductionn}
Let $G$ be a finite group,   $N$ be a
normal subgroup of $\,G$ and $\theta\in \Irr(N)$. Let 
$G_{\theta}$ be the stabilizer of $\,\theta$ in 
$\,G$.  
Then  $\eta(\theta^G)=\eta(\theta^{G_{\theta}})$.
\end{lem}
\begin{proof}
Observe that all the irreducible constituents of $\theta^{G_{\theta}}$ lie above 
$\theta$. Thus by Clifford theory it follows that $\eta(\theta^G)=\eta(\theta^{G_{\theta}})$.
\end{proof}
\begin{lem}\label{normal}
 Let $G$ be a finite $p$-group, $H$ be a normal 
 subgroup of $G$ and $\theta\in \Irr(H)$.
Then either $\eta(\theta^G)=1$ or $\eta(\theta^G)\geq p$.
\end{lem}
\begin{proof}
In Lemma 4.1 of \cite{edith}, it is proved that, if in addition of the previous 
hypotheses, $\theta$ is $G$-invariant, then $\eta(\theta^G)=1$ or $\eta(\theta^G)\geq p$.
Thus by induction on $|G:H|$ and Lemma \ref{inductionn}, the result follows.
\end{proof} 
Let $G$ be a group, $H$ be a subgroup of $G$ and $\theta\in \Irr(H)$. Denote
by
$\Irr(G\mid \theta)=\{\chi \in \Irr(G)\mid [\chi_H,\theta]\neq 0\}$ the set
or irreducible characters of $G$ lying above $\theta$.

\begin{lem}\label{morethanp}
Let $G$ be a finite $p$-group, $H$ be a subgroup of $G$ and $\theta \in \Irr(H)$. 
Let $Z_1$ be a subgroup of the center ${\bf Z}(G)$ such that $|HZ_1:H|=p$.
Then $\theta$ extends to $HZ_1$ and
\begin{equation*}
\eta(\theta^G)=\sum_{\nu\in \Irr(HZ_1\mid \theta)} \eta(\nu^G).
\end{equation*}
In particular, if  $\nu\in \Irr(HZ_1\mid \theta)$ we have that
\begin{equation}\label{extension}
\eta(\theta^G)\geq \eta(\nu^G)+ (p-1).
\end{equation}
\end{lem}
\begin{proof}
Observe that $\theta$ extends to $HZ_1$ since $Z_1\leq {\bf Z}(G)$ and $|HZ_1:H|=p$. 
Thus there are exactly $p$ characters in 
$\Irr(HZ_1\mid \theta)$. Let 
$\alpha\in \Lin(H\cap Z_1)$ be the unique character such that 
$\theta_{H\cap Z_1}= \theta(1) \alpha$. 
    Since
 $(\theta^{HZ_1})_{Z_1}=(\theta_{H\cap Z_1})^{Z_1}$, we have that
 $(\theta^{HZ_1})_{Z_1}=\theta(1) \sum_{\nu\in \Lin(Z_1\mid \alpha)}\nu$. Therefore
\begin{equation}\label{extensionlambda}
    \mbox{ for any }\nu,\mu \in \Irr(HZ_1\mid \theta)\mbox{, if } \nu\neq \mu \mbox{ then } \nu_{Z_1}\neq \mu_{Z_1}.
  \end{equation}  
   
    Observe that for any $\chi\in\Irr(G)$ and any $\beta \in \Lin(Z_1)$,
  if $[\chi_{Z_1},\beta]\neq 0$ then 
 $\chi_{Z_1}=\chi(1)\beta$. By \eqref{extensionlambda}, it follows that 
  if $\chi, \psi\in\Irr(G)$, $ \nu,\mu \in \Irr(HZ_1\mid \theta)$, $\nu\neq\mu$, 
  $[\chi_{Z_1},\nu]\neq 0$ and  $[\psi_{Z_1},\mu]\neq 0$, then 
 $\chi\neq \psi$. 
 Thus  the irreducible constituents of $\theta^G$ lying
over distinct extensions of $\theta$ in $HZ_1$ are distinct characters. 
It follows that 
\begin{equation*}
\eta(\theta^G)=\sum_{\nu\in \Irr(HZ_1\mid \theta)} \eta(\nu^G).
\end{equation*}
Since $\eta(\nu^G)\geq 1$ for any $\nu\in \Irr(HZ_1)$, \eqref{extension} follows.
\end{proof}
\end{section}
\begin{section}{Proof of Theorem A}

Let $G$ and $\theta\in \Irr(H)$ be a minimal counterexample of the statement of Theorem A
with respect to the order  $|G|$ of $G$. That is 
we are assuming that 
\begin{equation}\label{suposicion}
|G:H|=p^2,\, 1<\eta(\theta^G)<\frac{p+1}{2}
\end{equation}
\noindent and 
\begin{induction}\label{induc}
 for any  finite $p$-group $G_1$, any  subgroup $H_1$ of $G_1$, and any $\theta_1\in \Irr(H_1)$,
 if $\,|G_1:H_1|=p^2$ and $|G_1|<|G|$
  then either $\eta({\theta_1}^{G_1})=1$ or $\eta({\theta_1}^{G_1})\geq \frac{p+1}{2}$.
\end{induction}

 Set $\overline{L}=L/\core_G(\Ker(\theta))$ for any subgroup $L$ of $\,G$ such that
 $L\geq \core_G(\Ker(\theta))$. 
Observe that $H\geq \core_G(\Ker(\theta))$ and $|\overline{G}:\overline{H}|=|G:H|$.
 Observe also that we  can regard
$\theta$ as a character of $H/\core_G(\Ker(\theta))$ and $\eta(\theta^{\overline{G}})=\eta(\theta^G)$.

By working with the group $G/\core_G(\Ker(\theta))$ and \ref{induc},  
we may assume that $$\core_G(\Ker(\theta))=1.$$ 
Thus $\overline{L}=L$ for all subgroups $L$ of $G$.

Denote by $Z$ the center ${\bf Z}(G)$ of $G$. Let $\nu\in \Lin(Z)$ be the unique
character of $Z$ lying below $\theta$.

\begin{claim}\label{center}
$Z<H$. Thus $\nu\in \Lin(Z)$ is a multiple of a 
faithful character of $Z$ and 
$Z$ is a cyclic group.
\end{claim}
\begin{proof}
Suppose $Z$ is not contained in $H$. Let $Z_1\leq Z$ be such that $|HZ_1 :H|=p$.
Lemma \ref{morethanp} implies that $\eta(\theta^G)\geq p$, a contradiction with 
\eqref{suposicion}. Thus 
$Z\leq H$.

Since $\Ker(\theta)\cap Z$ is normal in $G$ and $\core_G(\Ker(\theta))=1$, it 
follows that $\theta_Z\in \Lin(Z)$ is a faithful character of $Z$. Therefore
$\nu\in \Lin(Z)$ is faithful and
$Z$ is cyclic.
\end{proof}

\begin{claim}\label{coreg}
$\core_G(H)=Z$. 
\end{claim}
\begin{proof}
Assume that there exists a normal subgroup $N$ of $G$ such that $N\leq H$ and $N/Z$ is 
a chief factor of $G$. Fix  $\beta\in \Irr(N)$ such that $[\theta_N, \beta]\neq 0$.   
Since $\nu \in\Lin(Z)$ is a faithful character, we can check that
${\bf C}_G(N)$ is a normal subgroup of $G$ of index $p$. Also the stabilizer $G_{\beta}$
of $\beta$ in $G$ is ${\bf C}_G(N)$. 

If $H\cap{\bf C}_G(N)< H$, by Clifford theory we have that there 
exists some $\alpha\in \Irr(H\cap {\bf C}_G(N))$ such that $\alpha^H=\theta$. 
Thus $\eta(\theta^G)=\eta(\alpha^G)$.
Since $|{\bf C}_G(N)|<|G|$ and $|{\bf C}_G(N):H\cap{\bf C}_G(N)|=p^2$, by \ref{induc} we have that
$\eta(\alpha^{{\bf C}_G(N)})=1$ or $\eta(\alpha^{{\bf C}_G(N)})\geq \frac{p+1}{2}$. By Lemma \ref{inductionn} we have then
that $\eta(\alpha^G)=1$ or $\eta(\alpha^G)\geq \frac{p+1}{2}$ and therefore
$\eta(\theta^G)=1$ or $\eta(\theta^G)\geq \frac{p+1}{2}$, a contradiction with \eqref{suposicion}.
 We may assume then
that $H<{\bf C}_G(N)$.

 Since $|{\bf C}_G(N):H|=p$,  
 $H$ is normal in ${\bf C}_G(N)$ and thus by Lemma \ref{normal} we have that either
 $\eta(\theta^{{\bf C}_G(N)})=1$ or $\eta(\theta^{{\bf C}_G(N)})=p$.
 By Lemma \ref{inductionn} and the previous statement, we have that 
$\eta(\theta^{G})=1$ or $\eta(\theta^G)\geq p$, a contradiction with \eqref{suposicion}.
 Thus 
 such $N$ can not exist and so $\core_G(H)=Z$.
\end{proof}
Let $Y/Z$ be a chief factor of $G$. 
By the previous claim, it follows that $HY>H$. Since $Y/Z$ has order $p$, 
we have that $|HY:H|=p$. Since $|G:H|=p^2$, it follows that $|G:HY|=p$ and thus
$HY$ is a normal subgroup of $G$. 

Let $C$ be the centralizer ${\bf C}_G(Y)=\{g\in G\mid y^g=y \mbox{ for all } y\in Y\}$ of $Y$ in $G$.

\begin{claim}\label{cindexp}
 $|G:C|=p$. Also, given 
any $\mu\in\Lin(Y)$ extension of the faithful character $\nu\in \Lin(Z)$,
 we have that 
the stabilizer $G_{\mu}$ of $\mu$ in $G$ is $C$.
\end{claim}
\begin{proof}
Since $\nu\in \Lin(Z)$ is a faithful character of the center of $G$ and
 $Y/Z$ is 
a chief factor of the $p$-group $G$, it follows that the index of the centralizer $C$
of $Y$ in $G$ is $p$. 
\end{proof}
\begin{claim}   $HY/Z$ is an elementary abelian 
$p$-group. Also, we may assume that ${\bf Z}(HY)\geq Y$ and thus 
$C=HY$. 
\end{claim}
\begin{proof}
Since  $|HY:H|=p$, we have that $(HY)'=\left\langle  [h,k]\mid h,k\in HY\right\rangle
\leq H$. Observe that $(HY)'$ is normal in $G$ since
$HY$ is normal in $G$ and $(HY)'$  is a characteristic subgroup of $HY$. 
Since $\core_G(H)=Z$, it
follows then that $(HY)'\leq Z$. Also, since
$Y/Z$ is of order $p$ and $Z<H$, $(HY)^p=\left\langle   k^p\mid k\in HY \right\rangle$ is a characteristic subgroup of 
the normal subgroup $HY$ of $G$ and it is
contained in $H$. It follows then that $(HY)^p\leq Z$ and thus $HY/Z$ is an elementary 
abelian $p$-group.

Observe that the center ${\bf Z}(HY)$ of $HY$ contains $Z$.
If $\,{\bf Z}(HY)=Z$, then there is a unique character in $\Irr(H)$ lying above
$\nu$ since $HY/Z$ is elementary abelian $p$-group and $\nu\in \Lin(Z)$ is a 
faithful character, and so $\eta(\theta^G)=1$ or $\eta(\theta^G)=p$, that is a contradiction with
\eqref{suposicion} and therefore it must follow that ${\bf Z}(HY)>Z$.
 By replacing $Y$ for a normal subgroup of $G$ contained in ${\bf Z}(HY)$
if necessary, we may assume then that $Y\leq {\bf Z}(HY)$ and thus ${\bf C}_G(Y)=HY$.
\end{proof}

 \begin{claim}\label{extendtohy}
The character $\theta\in \Irr(H)$ extends to $HY=C$. Thus $\theta^{C}$ is the 
sum of the $p$ distinct  extensions of $\,\theta$. 
\end{claim}
\begin{proof}
Since $|HY:H|=p$, we have that 
either $\theta^{HY}\in \Irr(HY)$
or $\theta^{HY}$ is the sum of the $p$ distinct   extensions of $\theta$. 

Suppose that $\theta^{C}\in \Irr(C)$. 
Let $\mu\in\Lin(Y)$ be the unique character of $Y$
 such that $[(\theta^{HY})_Y,\mu]\neq 0$. Since  
$G_{\mu}=C$, then $\theta^G\in \Irr(G)$. Thus 
$\theta^{HY}$ is the sum of the $p$ distinct   extensions of $\theta$. 
\end{proof}
Let $\rho_1,\ldots, \rho_p\in \Irr(HY)$ be the $p$ distinct 
 extensions of $\theta$. Since $|G:HY|=p$,  by Lemma  \ref{normal} 
  we must have
  that 
  \begin{equation}\label{etalambda}
  \rho_i^G\in \Irr(G).
  \end{equation} 
 Since $Z(C)\geq Y$, there is a unique character
 $\mu_i\in \Lin(Y)$ lying below 
 $\rho_i$. 
  
\begin{claim}\label{centercy} ${\bf Z}(C)=Y$.
\end{claim}
\begin{proof}
Clearly $Y\leq {\bf Z}(C)$.
Assume that $Y< {\bf Z}(C)$.     
 Let
$ X\leq {\bf Z}(C)$ 
  such that $X/Y$ is a chief 
factor of $G$ and  $Y< X\leq HY=C$. Observe that such $X$ exists since $HY$ is normal in 
$G$, and $X$ is abelian since $X\leq {\bf Z}(C)$.
We are going to conclude that $\nu\in \Lin(Z)$ is not a faithful character, which is a 
  contradiction with Claim
\ref{center}.

\begin{subclaim}\label{xgandy}
The subgroup $[X,G]$ generates $Y = [X, G]Z$ modulo $Z$. 
\end{subclaim} 
\begin{proof}
Since $ Y$ and $X$ are normal subgroups of $G$ with 
$Y\triangleleft X$ and $|X/Y|=p$, the chief factor $X/Y$ of the 
$p$-group $G$ is centralized by $G$. 
So $[X,G]\leq Y$. Suppose that $[X,G]Z< Y$.
Since $|Y/Z|=p$, we must have       
$[X,G]\leq Z= {\bf Z}(G)$.
So commutation  in $G$ induces a bilinear map 
$$d: x Z, g {\bf C}_G(X) \mapsto [x,g]$$
\noindent of $(X/Z)\times (G/ {\bf C}_G(X))$ into 
the cyclic group $Z$. This map $d$ is non-singular
on the right by the definition of ${\bf C}_G(X)$. It is non-singular
on the left since $Z={\bf Z}(G)$. 
 Because $\lambda \in \Lin( X \mid \nu)$
extends the faithful character $\nu\in \Irr(Z)$, and $|X:Z|=p^2$,
this implies that ${\bf C}_G(X)=G_{\lambda}$ has index $p^2$ in 
$G$. But $C$ fixes $\lambda$ since $X\leq {\bf Z}(C)$.
 Therefore $\lambda$ has, at the same 
time, $p^2$ distinct $G$ conjugates, and at most 
$p=|G:C|$ such conjugates. This contradiction proves the claim.
\end{proof}

Given any character $\rho\in\Irr(C)$, since $X\leq{\bf Z}(C)$, we have that
$\frac{1}{\rho(1)} \rho_X\in \Lin(X)$ is the unique character
lying below $\rho$.

\begin{subclaim} \label{3distinct}
 There exist some $\lambda\in \Lin(X)$,  
 some $g\in G\setminus C$ and $i\in \{2,\ldots,p-1\}$ such that $[(\theta^{C})_X,\lambda]\neq 0$,
 $[(\theta^{C})_X,\lambda^g]\neq 0$ and $[(\theta^{C})_X, \lambda^{g^i}]\neq 0$. 
\end{subclaim}
\begin{proof}
Since $1<\eta(\theta^G)< \frac{p+1}{2}$ and $\rho_1^G,\ldots, \rho_p^G$ are the
irreducible constituents of $\theta^G$,
 there exist at least 3 distinct 
  $j,k,l \in \{1,2,\ldots, p\}$ such that $\rho_j^G=\rho_k^G= \rho_l^G$.
  Since  $X$ is normal in 
  $G$, by Clifford Theory it follows that
  $\frac{1}{\rho_j(1)}(\rho_j)_X$, $\frac{1}{\rho_k(1)}(\rho_k)_X$ and 
  $\frac{1}{\rho_l(1)}(\rho_l)_X$
    are $G$-conjugates.  Set $\lambda=\frac{1}{\rho_j}(\rho_j)_X$. 
 Then there
  exists some $g\in G\setminus C$ 
  such that $\lambda^g=\frac{1}{\rho_k(1)}(\rho_k)_X$
   Since $X\leq {\bf Z}(C)$ and $|G:C|=p$,
  there exists some $i\in \{2,\ldots, p-1\}$ such that
  $(\lambda)^{g^i}=\frac{1}{\rho_l(1)}(\rho_l)_X$.
 \end{proof}
 
Fix $g\in G\setminus C$  as in  \ref{3distinct}.
Since $X/Y$ is cyclic of order $p$, $H\cap X> Z$, and $H\cap Y=Z$ we may choose 
\begin{equation}
x\in H \mbox{ such that } X=\left\langle x,Y\right\rangle
\end{equation}
   Since $X\leq {\bf Z}(C)$, we have $[X,C]=1$. 
Suppose that $[x,g^{-1}]\in Z$. Then $x$ centralizes both 
$g^{-1}$ and $C$ modulo $Z$. Hence $xZ \in {\bf Z}(G/Z)$,
which is false by Step \ref{xgandy}. Hence 
$[x, g^{-1}] \in Y\setminus Z$ and so  
\begin{equation}\label{18}
Y=Z\left\langle   y \right\rangle  \mbox{ is generated over } Z \mbox{ by } y=[x,g^{-1}].
\end{equation}
Since $[Y,G]\leq Z$ we have that $z=[y,g^{-1}]  \in Z$. If
$z=1$, then $G=C\left\langle g\right\rangle$ centralizes $Y = Z\left\langle y\right\rangle$, since $C$ centralizes 
$Y<X$ because  $X\leq {\bf Z}(C)$, and $G$ centralizes Z. This
is impossible because $Z={\bf Z}(G)< Y$. Thus 
\begin{equation}\label{19}
z= [y,g^{-1}] \mbox{ is a non-trivial element of }Z.
\end{equation}
By \eqref{18} we have  $y=[x,g^{-1}]=x^{-1} x^{g^{-1}}$. By
\eqref{19} we have $z= [y, g^{-1}]=y^{-1} y^{g^{-1}}$ . Finally
$z^{g^{-1}}= z$ since 
$z\in Z$. Since $X=Z\left\langle x,y  \right\rangle$ is abelian since $X\leq {\bf Z}(C)$,
it follows that 
\begin{equation}\label{21}
z^{g^{-j}}=z, \  y^{g^{-j}}= yz^j \mbox{ and }
 x^{g^{-j}}=xy^j z^{\binom{j}{2}},
\end{equation}
\noindent for any integer $j= 0, 1, \ldots, p-1$. 
Because $g^{-p} \in C$ centralizes $X$ since $X\leq {\bf Z}(C)$, 
we have
\begin{equation*} 
z^p=1 \mbox{ and } y^p z^{\binom{p}{2}}=1.
\end{equation*}
Since $p>2$ is odd by hypothesis, $p$ divides
$\binom{p}{2}= \frac{p(p-1)}{2}$ and $z^{\binom{p}{2}}=1$. Therefore
$y^p= z^p =1$.
It follows that $y^i$, $z^i$ and $z^{\binom{i}{2}}$ 
depend only on the residue of $i$ modulo 
$p$, for any integer $i\geq 0$.
such 
that $X = Y    \left\langle  x \right\rangle$ and $x\in C$. Thus by \eqref{19}
we have that 
\begin{equation}\label{not1} 
z^{\binom{j}{2}}\neq 1 \mbox{ for any integer } 0<j<p.
\end{equation}

Let  $\lambda\in \Lin(X)$ and $i\in\{2,\ldots, p-1\}$ be
as in Step
 \ref{3distinct}. Set  $\varpi= \frac{1}{\theta(1)}\theta_{X\cap H}$.
 We can check that $\varpi\in \Lin (X\cap H)$. Since $(\theta^C)_X=(\theta_{H\cap X})^X$,
 we have that 
 $\lambda$, $\lambda^g$ and $\lambda^{g^i}$ are extensions of
 $\varpi$. Since  $x\in (H\cap X)$, by the previous statement
 we have that
 \begin{equation}\label{equal}
 \lambda(x)= \lambda^g(x)=\lambda^{g^i}(x).
 \end{equation}
By \eqref{21} we have that
  $$\lambda^g(x)= \lambda(x^{g^{-1}})=\lambda(xy)=\lambda(x)\lambda(y).$$
  Thus by \eqref{equal}, we get
  \begin{equation}\label{y1}
  \lambda(y)=1.
  \end{equation}
  
 Therefore
 \begin{eqnarray*}
   \lambda^{g^{i}}(x)& = & \lambda(x^{g^{-i}})\\
 &= &  \lambda(xy^{i}z^{\binom{i}{2}}) \ \ \mbox{ by \eqref{21}}\\
  &=&  \lambda (x)\lambda(y^{i})\lambda(z^{\binom{i}{2}})\\
  & =& \lambda(x) \lambda (z^{\binom{i}{2}}),      
  \end{eqnarray*}
  \noindent where the last line follows from \eqref{y1}. By \eqref{equal}, we have that
   $\lambda (z^{\binom{i}{2}})=1$. But $\lambda_Z=\nu\in \Lin(Z)$ 
   is a faithful character and 
   $z^{\binom{i}{2}}\neq 1$ by \eqref{not1}.
    This is a contradiction and the claim is proved. 
 \end{proof}
 Since ${\bf Z}(HY)=Y$, we have that ${\bf Z}(H)=Z$. Thus $HY$ is a class 2 group with 
 $HY/Z$ elementary abelian. Therefore $\theta\in \Irr(H)$ is the only character in $H$ lying above
 $\nu\in \Lin(Z)$. Hence an irreducible character of $G$ lies over $\theta$ if and 
 only if it lies over $\nu$. Since $\Irr(G\mid \nu)$ has either 1 element or
 at least $p$ by Lemma \ref{normal}, 
 it follows that $\eta(\nu^G)=1$ or $\eta(\nu^G)\geq p$, and therefore
 either $\eta(\theta^G)=1$ or 
 $\eta(\theta^G)\geq p$. But $1<\eta(\theta^G)<\frac{p+1}{2}$, and that
 is our final contradiction and thus the statement of Theorem A holds.  
\end{section}
\begin{section}{Examples}

In this section, we will prove that the group $G$, the subgroup $H$ and
the character 
$\lambda\in \Lin(H)$ that satisfy Hypotheses \ref{dadeexample}
have the properties that $|G:H|=p^2$ and $\eta(\lambda^G)=\frac{p+1}{2}$. 
And then, given any integer $n\geq 2$, we  construct a group $G$ with a subgroup 
$H$ and a character $\lambda\in \Lin(H)$ such that $|G:H|=p^n$ and $\eta(\lambda^G)=\frac{p+1}{2}$.
 
 \begin{hyp}\label{dadeexample} Fix an odd prime $p$.
Let $G$ be the semidirect product of a cyclic group $C$ of order $p$ 
and an elementary abelian 
group  $A$ of order $p^3$. Assume $C=\left\langle  c \right\rangle$ and 
\begin{equation}
A=\left\langle  a \right\rangle\times \left\langle [a,c]\right\rangle\times \left\langle[a,c,c]\right\rangle,
\end{equation}
\noindent for some $a$ in $A$.
Observe that the subgroup $\{e\}\times\{e\}\times \left\langle[a,c,c]\right\rangle$
 is the center of the group $G$.
Set  $Z=\{e\}\times\{e\}\times \left\langle [a,c,c]\right\rangle$. 
 
Fix  $\omega$   a primitive complex $p$-th root of unity. 
Let $\alpha\in\Lin(\left\langle a\right\rangle)$, $\beta\in\Lin( \left\langle 
[a,c]\right\rangle)$ and 
$\gamma\in \Lin(\left\langle [a,c,c]\right\rangle)$ be the unique linear characters such that
$\alpha(a)=\beta([a,c])=\gamma([a,c,c])=\omega$.

Set 
\begin{equation}
H=\left\langle  a   \right\rangle\times \{ e\} \times \left\langle  [a,c,c]\right\rangle\mbox{ and } \lambda=1_{\left\langle  a\right\rangle}\times 1_{\{e\}}\times \gamma\in \Lin(H).
\end{equation}

Observe that $H$ is a subgroup of $A$ of index $p$. Thus $|G:H|=p^2$.
Observe also that $\lambda$ extends to $A$ and there are exactly $p$ distinct
extensions of $\lambda$ to $A$, namely
\begin{equation}
\Irr(A\mid \lambda)=\{1_{\left\langle  a \right\rangle}\times \beta^r \times \gamma\mid r=0,1,\ldots, p-1\}.
\end{equation}
Set $\Lambda_{r}= 1_{\left\langle a\right\rangle}\times \beta^r \times \gamma$.
\end{hyp}
\begin{lem}\label{counting} Assume Hypotheses \ref{dadeexample}
 Given any integer $i$ with $0< i$, we have that 
$$(\Lambda_{r})^{c^i}= \alpha^{ri+ \frac{i (i-1)}{2}}\times \beta^{r+i}\times\gamma.$$
\end{lem}
\begin{proof}
Observe that $(\Lambda_{r})^c=\alpha^r\times\beta^r\beta\times\gamma=\alpha^r\times \beta^{r+1}\times \gamma$ since
$a^c=a[a,c]$ and $[a,c]^c=[a,c][a,c,c]$. 
Assume by induction that  
\begin{equation}\label{countingn}
(\Lambda_{r})^{c^n}= \alpha^{rn+ \frac{n (n-1)}{2}}\times \beta^{r+n}\times\gamma.
\end{equation}
Then 
\begin{eqnarray*}
(\Lambda_{r})^{c^{n+1}}&= &((\Lambda_{r})^{c^n})^c\\
&=& (\alpha^{rn+ \frac{n (n-1)}{2}}\times \beta^{r+n}\times\gamma)^c \mbox{ by \eqref{countingn} }\\
&=& \alpha^{rn+ \frac{n (n-1)}{2}+ r+n }\times \beta^{r+n+1}\times\gamma,
\end{eqnarray*}
\noindent where the last line follows since
$a^c=a[a,c]$ and $[a,c]^c=[a,c][a,c,c]$. We can check
that $rn+ \frac{n (n-1)}{2}+ r+n= r(n+1) + \frac{(n+1)(n)}{2}$. Thus 
\begin{equation*}
(\Lambda_{r})^{c^{n+1}}=\alpha^{r(n+1)+ \frac{(n+1)n}{2}}\times \beta^{r+(n+1)}\times\gamma,
\end{equation*} 
\noindent and the result follows by induction.
\end{proof}
\begin{lem}\label{randi} Assume Hypotheses \ref{dadeexample}.
Let $r$ be an integer such that $0<r<p$. 
Then $(\Lambda_r)^{c^j}$ is an extension of
$\lambda$ if and only if either $j\equiv 0\modu p$ or $j\equiv (1-2r) \modu p$. 
If   $\,i\equiv (1-2r) \modu p$ then  $(\Lambda_r)^{c^i}= \Lambda_{1-r}$. 
\end{lem}

\begin{proof}
By Lemma \ref{counting}, we have that
  $(\Lambda_r)^{c^i}$ is an extension of $\lambda$ if and only if 
$\alpha^{ir+\frac{i(i-1)}{2}}=1_{\left\langle  a\right\rangle}$. 
 Since $\alpha$ is a faithful linear character of a cyclic group of order $p$,   $\alpha^{ir+\frac{i(i-1)}{2}}=1_{\left\langle a\right\rangle}$ if and only if
   $(ir+\frac{i(i-1)}{2})\equiv 0\modu p$.
Observe that  $(ir+\frac{i(i-1)}{2})\equiv 0\modu p$ if and only if either $i\equiv 0\modu p$ or
$(r+\frac{i-1}{2})\equiv 0 \modu p$. Therefore
$(\Lambda_r)^{c^i}$ is an extension of $\lambda$  if and only if either $i\equiv 0\modu p$ or
$i\equiv(1-2r) \modu p$. 

If $i\equiv(1-2r)\modu p$, then $(\Lambda_r)^{c^i}= \Lambda_{1-r}$  by Lemma \ref{counting}.
\end{proof}
 
\begin{lem}\label{inbetween} Assume Hypotheses \ref{dadeexample}.
Then 
 $1<\eta(\lambda^G)\leq\frac{p+1}{2}.$
\end{lem}
\begin{proof}
By the previous lemma, it follows that the stabilizer of  
$\Lambda_{r}$ is a proper subgroup of
$G$. Since $|G:A|=p$ and $\Lambda_r\in\Lin(A)$, we have that 
\begin{equation}
  (\Lambda_{r})^G\in \Irr(G) \mbox{ for any integer } r.
 \end{equation}   
 
Since $p<2$, it follows that there exists two distinct integers $k,l$ such that
 $0<k, l<p$ and $k\neq (1-2l)\modu p$. Thus by Lemma \ref{randi}
 we have that $\Lambda_k$ and $\Lambda_l$ are
 not $G$-conjugates. It follows that $(\Lambda_k)^G\neq(\Lambda_l)^G$.
 Since $(\Lambda_k)^G\neq(\Lambda_l)^G$, $(\Lambda_k)^G,(\Lambda_l)^G\in \Irr(G)$
 and  both
 $\Lambda_k$ and $\Lambda_l$ lie above $\lambda$,
 we have that $\eta(\lambda^G)\geq 2$.
 
  Observe that $r\equiv(1-r) \modu p$  if and only ir $2r\equiv1\modu p$. Thus given any
  $r$ such that $0<r<p$ and $2r\neq 1\modu p$, by Claim \ref{randi} we have that
  $\Lambda_{r},\Lambda_{1-r}\in \Irr(A)$ are two distinct 
  $G$-conjugate extensions of $\lambda$. Thus $\eta(\lambda^G)\leq \frac{p+1}{2}$.
 \end{proof}

\begin{prop}\label{dade}
Assume Hypotheses \ref{dadeexample}.
Then $|G:H|=p^2$ and $\eta(\lambda^G)=\frac{p+1}{2}$.
\end{prop}
\begin{proof}
By  Lemma \ref{inbetween},  we have that $1<\eta(\lambda^G)\leq \frac{p+1}{2}$.
Thus by Theorem A, it follows that $\eta(\lambda^G)=\frac{p+1}{2}$.
\end{proof}
Denote by $1_H$ the principal character of $H$. 
\begin{lem}\label{dade2}
Let $p$ be a prime number, $G$ be a $p$-group and $H$ be a subgroup of $G$ with 
$|G:H|=p^n$.
Then $\eta((1_H)^G)\geq n(p-1)+1$.
\end{lem}
\begin{proof}
Using induction on the order of $G$, without lost of generality we may assume
that $\core_G(H)=1$. We are going to use induction on $n$.

Let $Z_1$ be a subgroup of the center ${\bf Z}(G)$
of $G$ with $|Z_1|=p$. Observe that $H\cap Z_1=1$ since $\core_G(H)=1$.
Thus $|HZ_1:H|=p$.
By Lemma \ref{morethanp}, we have that
\begin{equation}\label{more}
\eta((1_H)^G)\geq \eta((1_{HZ_1})^G)+(p-1).
\end{equation}
 Since $|G:HZ_1|=p^{n-1}$, by induction on  $n$ we have that
$$\eta( (1_{HZ_1})^G)\geq(n-1)(p-1)+1.$$
 The result follows by \eqref{more} and the previous statement.
\end{proof}

\begin{lem}\label{basico}
Let $G_0$ be a $p$-group and $\Gamma$ be a character of $G_0$. Assume that 
$[\Gamma, 1_{G_0}]=0$.
Let $N=G_0\times G_0\times \cdots\times G_0$ be the direct product of $p$-copies of $G_0$.
Set $$\Delta=\Gamma\times 1_{G_0}\times\cdots\times 1_{G_0}.$$ Let
$C=\left\langle c \right\rangle$ be a cyclic group of order $p$. 
 Observe that $C$ acts on $N$ by
 \begin{equation}\label{actionc}
 c:(n_0, n_1,\ldots,n_{p-1})\mapsto(n_{p-1},n_0, \dots, n_{p-2})
 \end{equation}
 \noindent for any $(n_0, n_1, \ldots,n_{p-1})\in N$.
 
Let $G$ be the direct product of $N$ and $C$, i.e $G$ is the wreath product of $G_0$ and 
$C$.  Then $\eta(\Delta^G)=\eta(\Gamma)$. 
\end{lem}
\begin{proof}
Let $\delta\in\Irr(N)$ be a constituent of $\Delta$. Observe that
$\delta$ is of the form $\gamma\times 1_{G_0}\times\cdots\times 1_{G_0}$, for
some $\gamma\in \Irr(G_0)$ such that $[\gamma,\Gamma]\neq 0$. 
Observe that $\gamma\neq 1_{G_0}$ since $[\Gamma, 1_{G_0}]=0$.
By \eqref{actionc}, we have that
$\delta$ is $G$-invariant if and
only if $\gamma=1_{G_0}$.
Thus  $\delta^G\in \Irr(G)$ for any constituent $\delta\in \Irr(N)$ of $\Delta$.
 Observe that the $G$-orbit of $\delta\in \Irr(N)$ is 
 $$\{\gamma\times 1_{G_0}\times\cdots \times 1_{G_0}, 1_{G_0}\times\gamma\times\cdots\times 1_{G_0},\cdots, 1_{G_0}\times\ldots\times 1_{G_0}\times \gamma\}.$$
Thus if $ \delta, \epsilon \in \Irr(N)$ are two distinct constituents of $\Delta$, then
 $\delta^G\neq \epsilon^G$. It follows that $\eta(\Delta^G)=\eta(\Gamma)$.
\end{proof}
\begin{theo}\label{extensiondade}
Let $p$ be an odd prime number and $n\geq 2$ be an integer. There exist a $p$-group
$G$, a subgroup $H$ of $\,G$ and $\lambda\in \Lin(H)$,
such that $|G:H|=p^n$
  and  $\eta(\lambda^G)=\frac{p+1}{2}$.
\end{theo}
\begin{proof}
 If $n=2$, then the result follows by Lemma \ref{dade}.
  By induction on $n$, we may 
 assume that the result holds for any integer $n$ such that  $n-1\geq 2$.
 \begin{com}\label{induction}
  Fix a $p$-group $G_0$, a subgroup $H_0\leq G_0$
 and  
 $\lambda_0\in \Lin(H_0)$ such that $|G_0:H_0|=p^{n-1}$
 and  $\eta(\lambda_0^{G_0})=\frac{p+1}{2}$.
 \end{com}

 Let $N$ and $G$ be as in Lemma \ref{basico}.
 Let
 $$H=H_0\times G_0\times\ldots\times G_0.$$ Then $H$ is a subgroup of 
 $N$ and $|G:H|=|G:N||N:H_0|=p|G_0:H_0|=p^n$.  
 
 Set $\lambda=\lambda_0\times 1_{G_0}\times \ldots\times 1_{G_0}$. Observe that
 $\lambda\in \Lin(H)$ since $\lambda_0\in \Lin(H_0)$.
  We can check that 
 $\eta(\lambda^N)=\eta(\lambda_0^{G_0})$. Thus by \ref{induction} we have that
 $\eta(\lambda^N)=\frac{p+1}{2}$. 

By Lemma \ref{dade2}, we have that $\lambda_0\neq 1_{H_0}$. 
Thus $[\lambda_0^{G_0}, 1_{G_0}]=0$.
By Lemma \ref{basico} we have then that 
 $\eta(\lambda^N)=\eta(\lambda^G)$
and the result is proved.
\end{proof}

\begin{lem}\label{pcube}
Let $p$ be a prime number such that $p-1$ is divisible by 3.
Fix $r\in \{1,\ldots, p-1\}$. 
Then the set  $\{r(1- i^3)\modu p \mid i=0,\ldots ,p-1\}$ has  
$\frac{p+2}{3}$ elements. Also, given any $e\in \{r(1- i^3)\modu p \mid i=1,\ldots ,p-1\}$,
 there are exactly 3 distinct solutions in  $\{1,\ldots, p-1\}$ of
the equation $e\equiv r(1-x^3)\modu p$
\end{lem}
\begin{proof} Let $u$ be a generator of the units of the field $Z_p$ of 
$p$ elements. Then $U=\left\langle  u^{\frac{p-1}{3}}\right\rangle$ is a subgroup of 
order 3 and any element in $U$ is a solution of $x^3\equiv 1 \modu p$. 
Thus given any
integer $n\neq r$, if the equation $x^3\equiv r-n\modu p$ has a solution, then it has
exactly 3 distinct solutions in $Z_p$. Therefore the set
$ \{r(1- i^3)\modu p \mid i=1,\ldots ,p-1\}$ has $\frac{p-1}{3}$ distinct elements.
Since $0^3=0$, the set $\{ (r(1-i^3)\modu p \mid i=0,\ldots ,p-1\}$  has 
$\frac{p-1}{3}+1=\frac{p+2}{3}$.
\end{proof}
\begin{hyp}\label{construction2}
Let $p>5$ be a prime number such that $p-1$ is divisible by $3$.
 Let $F$ be a field of $p$ elements and $F[x]$ be the truncated
polynomial algebra generated over $F$ by some $x$ satisfying only $x^4=0$. So 
$F[x]$ is a vector space of dimension 4 over $F$ with $1$, $x$, $x^2$ and $x^3$ as a
basis. Let $m$ be an isomorphism of the additive group $F[x]^{+}$ of $F[x]$ onto a 
multiplicative group $M$. Then $M$ is an elementary abelian multiplicative group of 
order $p^4$ with $m(1)$, $m(x)$, $m(x^2)$, $m(x^3)$ as generators. Let $U$ be the subgroup of
the unit group $F[x]^{\times}$ generated by $1+x$ and $1+x^2$. The general element of 
$U$ is 
\begin{equation}\label{one}
(1+x)^i (1+x^2)^j= 1+ix+(\binom{i}{2} +j) x^2 + (\binom{i}{3}+ij)x^3
\end{equation}
\noindent for arbitrary integers $i$, $j$, since $x^4=0$. Because $p>3$, it follows that
$U$ is elementary abelian of order $p^2$, and that \eqref{one} holds for any $i,j\in F$.
The group $U$ acts naturally on the group $M$, so that 
\begin{equation}\label{two}
m(y)^u=m(yu)
\end{equation}
\noindent for all $y\in F[x]$ and $u\in U$. Let $G$ be the semidirect product of $M$ and $U$.
Then $G$ is a multiplicative group with order $p^6$.   

Let $H$ be the subgroup 
\begin{equation}\label{defih}
H=\left\langle  m(1), m(x),m(x^3) \right\rangle=\{m(a_0+a_1x+a_3 x^3)\mid a_0, a_1, a_3 \in F\}.
\end{equation}

Fix a primitive $p$-th root of unity $\omega$. Fix an integer $r>0$ such that
$3r\equiv -1 \modu p$. Thus $r\equiv \frac{-1}{3}\modu p$  and $r\not\equiv 0\modu p$. 
Let $\lambda\in \Lin(H)$ be the 
character given by 
\begin{equation}\label{defilambda}
\lambda(m(a_0+a_1x+a_3 x^3)) =\omega^{ra_0+ra_1+a_3}.
\end{equation}

\end{hyp}

\begin{theo}\label{examples2}
Assume Hypothesis \ref{construction2}. Then 
\begin{equation}
\lambda^G= \chi_0 + 3\sum_{i=1}^{\frac{p-1}{3}} \chi_i
\end{equation}
 \noindent where $\chi_i\in \Irr(G)$ and  $\chi_i\neq \chi_j$ if 
 $i\neq j$ for $i=0,1, \ldots , \frac{p-1}{3}$. Thus 
$\eta(\lambda)=\frac{p+2}{3}$.
\end{theo}
\begin{proof}
The center ${\bf Z}(G)$ of $G$ is the subgroup $\left\langle m(x^3) \right\rangle$ or 
order $p$.  Let 
$\gamma$ be the faithful linear character of ${\bf Z}(G)$ sending $m(x^3)$ to $\omega$. 
Then $\Lin (M\mid \gamma)$ consisting of the $p^3$ linear characters
$\mu_{f_0, f_1, f_2}$, for $f_0,f_1, f_2\in F$ given by 
\begin{equation}
\mu_{f_0, f_1, f_2}( m(a_0+a_1x+a_2x^2+a_3x^3))= \omega^{f_0 a_0+ f_1a_1+f_2a_2+a_3}
\end{equation}
for all $a_0$, $a_1$, $a_2$, $a_3\in F$. If $e,i,j \in F$, then \eqref{one} and \eqref{two}
imply that the conjugate character $ \mu_{e, 0, 0}^{(1+x)^{-i}(1+x^2)^{-j}}$ to
$\mu_{e,0,0}$ sends

\begin{eqnarray*}
m(1)& \mapsto & \mu_{e,0,0}(m(1+ix+(\binom{i}{2} +j)x^2 + (\binom{i}{3}+ij)x^3))=\omega^{e+\binom{i}{3}+ij},\\
m(x)&\mapsto & \mu_{e,0,0}(m(x+ix^2 + (\binom{i}{2}+j)x^3))=\omega^{\binom{i}{2}+j},\\ 
m(x^2)&\mapsto & \mu_{e,0,0}(m(x^2 + ix^3))=\omega^i,\\ 
m(x^3) &\mapsto & \mu_{e,0,0}(m(x^3))=\omega.
\end{eqnarray*}
\noindent It follows that
\begin{equation}\label{cuatro}
 \mu_{e, 0, 0}^{(1+x)^{-i}(1+x^2)^{-j}}=\mu_{e+\binom{i}{3}+ij, \binom{i}{2}+j, i}
 \end{equation}
 for any $e,i,j\in F$. If we fix $e$, then the above equation implies that distinct
 pairs $(i,j)\in F\times F$ yield distinct conjugates
   $ \mu_{e, 0, 0}^{(1+x)^{-i}(1+x^2)^{-j}}\in \Lin(M\mid \gamma)$.
    Hence the $G$-orbit $L_e$ of
 $\mu_{e,0,0}$ has exactly $p^2$ members. 
 Furthermore the above equation implies that the only member of that orbit with the form 
 $\mu_{f,0,0}$ is $\mu_{e,0,0}$. We conclude that the orbit $L_e$, for $e\in F$, 
 are $p$ distinct $G$-orbits in $\Lin(M\mid\gamma)$, each with size $p^2$.
 Since the normal subgroup $M$ of index $p^2$ is exactly the stabilizer of $\mu_{e,0,0}\in \Lin(M)$ in $G$, the induced characters 
 \begin{equation}\label{five}
 \chi_e=\mu_{e,0,0}^G \mbox{ are precisely the distinct members of }\Irr(G\mid \gamma).
 \end{equation} 
 \noindent Then 
 \begin{equation}\label{ocho}
 \lambda^M=\sum_{f\in F} \mu_{r,r,f} \mbox{ and } \lambda^G= \sum_{f\in F} \mu_{r,r,f}^G.
  \end{equation}
  \begin{claim}
  Let $i\in \{1,\ldots ,p-1\}$, $e= r(1-i^3)$ and $j= r- \binom{i}{2}$. Then 
  \begin{equation}\label{conjugate}
  \mu_{e,0,0}^{(1+x)^{-i}(1+x^2)^{-j}}= \mu_{r,r,i}. 
  \end{equation}
  \end{claim}
 \begin{proof}
  For a fix $i$, we have 
  \begin{eqnarray*}
   e+\binom{i}{3}+ij &=& e +\binom{i}{3} + i(r-\binom{i}{2})\\
   &=& e+ \frac{i(i-1)(i-2)}{6} +i(r-\frac{i(i-1)}{2})\\
   &=& i^3(\frac{1}{6}-\frac{1}{2}) +i^2(\frac{1}{2}-\frac{1}{2}) + i(r+\frac{1}{3})+e\\
   &\equiv & \frac{-i^3}{3} +e\modu p, \mbox{\,\,\, since $r\equiv \frac{-1}{3} \modu p$ }\\
   &\equiv  & \frac{-i^3}{3} + r(1-i^3) \modu p, \mbox{\,\, since $e= r(1-i^3)$}\\
   &\equiv &  r - i^3( r+\frac{1}{3})\equiv r\modu p,
   \end{eqnarray*}
   \noindent where the last line follow since $r\equiv \frac{-1}{3} \modu p$. Thus
    $(e+\binom{i}{3}+ ij, \binom{i}{2}+j, i)= (r,r,i)$ in $F\times F\times F$ and so by
    \eqref{cuatro} we get 
    \eqref{conjugate}. 
   \end{proof} 
   By the previous claim and \eqref{ocho}, we have that
    \begin{equation*}
  \lambda^G=\sum_{i=0}^{p-1} \mu_{r(1-i^3),0,0}^G.
  \end{equation*}
   By Lemma \ref{pcube}, we have then
\begin{equation}\label{nueve}
  \lambda^G= \mu_{r,0,0}^G+  3\sum_{e\in \{ r(1-i^3)\mid i=1,\ldots, p-1\}} \mu_{e,0,0}^G.
  \end{equation}
    By \eqref{five} we have that $\mu_{e,0,0}^G\in \Irr(G)$ and 
    $\mu_{e,0,0}^G\neq\mu_{f,0,0}^G$ if $e\not\equiv f \modu p$. Thus by Lemma \ref{pcube}
    and \eqref{nueve}, we conclude that 
    $\eta(\lambda^G)=\frac{p+2}{3}$ and the proof is complete.
\end{proof}
 \end{section} 
 {\bf Acknowledgment.} Professor Everett C. Dade brought to my attention
   Proposition \ref{dade}, the main step in the proof of Theorem \ref{extensiondade}, 
    and
   Lemma  \ref{dade2}. The proofs of those are based on his arguments,
    which are including here with his permission. 
   I thank him for that and for very useful conversations
   and emails.
   I also thank the referee of a previous version of this note for helpful comments and suggestions
   and for providing  Theorem \ref{examples2} when $p=7$, which is included here along with 
   the suggestions with
   her/his permission.

\end{document}